\documentclass[12pt,a4paper]{amsart}
\usepackage{amssymb,amsmath,amsthm}
\usepackage{graphicx}
\usepackage{hyperref}
 \usepackage[color=green!40]{todonotes}
\usepackage{mathtools}
\usepackage{enumerate}

\newcommand\cB{{\mathcal B}}
\newcommand\cC{{\mathcal C}}
\newcommand\cD{{\mathcal D}}

\newcommand\cF{{\mathcal F}}

\newcommand\cG{{\mathcal G}}

\newcommand\cU{{\mathcal U}}

\newcommand{\sat}{\mathop{}\!\mathrm{sat}}

\newcommand{\tops}{ \mathrel{\mathrlap{\top}\hspace{0.5ex}\top}
}
\newcommand{\bots}{ \mathrel{\mathrlap{\bot}\hspace{0.5ex}\bot}
}
\newcommand{\topbot}{\mathrel{\mathrlap{\tops}\hspace{-0.65ex}\bots}}

\newcommand\pisat{\mathop{}\!\scalebox{0.66}{$\tops$}\!\sat}

\newcommand\extsat{\mathop{}\!\scalebox{0.66}{$\topbot$}\!\sat}

\makeatletter
\newtheorem*{rep@theorem}{\rep@title}
\newcommand{\newreptheorem}[2]{%
\newenvironment{rep#1}[1]{%
 \def\rep@title{#2 \ref{##1}}%
 \begin{rep@theorem}}%
 {\end{rep@theorem}}}
\makeatother

\theoremstyle{plain}
\newtheorem{theorem}{Theorem}[section]
\newreptheorem{theorem}{Theorem}

\newtheorem{corollary}[theorem]{Corollary}
\newtheorem{conjecture}[theorem]{Conjecture}
\newtheorem{proposition}[theorem]{Proposition}

\theoremstyle{definition}

\newtheorem{defn}[theorem]{Definition}

\newtheorem{claim}[theorem]{Claim}

\newcommand\cref[1]{Corollary~\ref{cor:#1}}

\title{Projective and external saturation problem for posets}

\author{D\"om\"ot\"or P\'alv\"olgyi}
\address{ELTE Eötvös Loránd University and Alfréd Rényi Institute of Mathematics, Budapest}
\email{dom@cs.elte.hu}
\thanks{DP was partially supported by the ERC Advanced Grant ``ERMiD'' and by the J\'anos Bolyai Research Scholarship of the Hungarian Academy of Sciences, and by the New National Excellence Program \'UNKP-22-5 and by the Thematic Excellence Program TKP2021-NKTA-62 of the National Research, Development and Innovation Office.}

\author{Bal\'azs Patk\'os}
\address{Alfr\'ed R\'enyi Institute of Mathematics, Budapest}
\email{patkos@renyi.hu}
\thanks{Patk\'os's research is partially supported by NKFIH grants SNN 129364 and FK 132060.}

\date{}

\begin{document}
\begin{abstract}
We introduce two variants of the poset saturation problem. For a poset $P$ and the Boolean lattice $\cB_n$, a family $\cF$ of sets, not necessarily from $\cB_n$, is \textit{projective $P$-saturated} if (i) it does not contain any strong copies of $P$, (ii) for any $G\in \cB_n\setminus \cF$, the family $\cF\cup \{G\}$ contains a strong copy of $P$, and (iii) for any two different $F,F'\in\cF$ we have $F\cap[n]\neq F'\cap [n]$. Ordinary strongly $P$-saturated families, i.e., subfamilies $\cF$ required to be from $\cB_n$ satisfying (i) and (ii), automatically satisfy (iii) as they lie within $\cB_n$. We study what phenomena are valid both for the ordinary saturation number $\sat^*(n,P)$ and the projective saturation number $\pisat(n,P)$, the size of the smallest projective $P$-saturated family. 

Note that the projective saturation number might differ for a poset and its dual.
We also introduce an even more relaxed and symmetric version of poset saturation, \textit{external saturation}. We conjecture that all finite posets have bounded external saturation number, and prove this in some special cases.
\end{abstract}

\maketitle


\section{Introduction}

In this short note, we add several new notions and some related results to the growing literature of poset saturation problems. $\cB_n=(2^{[n]},\subseteq)$ denotes the Boolean poset of dimension $n$.

We say that a family $\cG\subseteq \cB_n$ is a \textit{weak copy} of a poset $P$ if there exists a bijective poset-homomorphism from $P$ to $\cG$, and $\cG$ is a \textit{strong copy} of $P$ if there exists a poset-isomorphism between $P$ and $\cG$. A family that does not contain a weak / strong copy of $P$ is called \textit{weak / strong $P$-free}. The extremal numbers $La(n,P)$ and $La^*(n,P)$ denote the maximum number of sets in a weak / strong $P$-free family $\cF\subseteq \cB_n$. In the past four decades, there have been a huge interest in determining (the asymptotics of) these parameters. For a survey, see \cite{GL} or \cite[Chapter 7]{GP}. 

The first instance of the corresponding saturation problem was introduced in \cite{Getal}, and the general problem was formalized in \cite{Fetal}: Let $\sat(n,P)$ / $\sat^*(n,P)$ denote the minimum size of a weak / strong $P$-free family $\cF\subseteq \cB_n$ such that $\{G\}\cup \cF$ contains a weak / strong copy of $P$ for any $G\in \cB_n\setminus \cF$. So far research has mainly focused on the cases when $P$ is an antichain \cite{Betal, DI, Fetal}, a chain \cite{Getal, MNS}, from a special subclass of posets \cite{Fetal}, or sporadic particular posets \cite{Fetal, I, I2, MSW}.
A general dichotomy phenomenon was stated first in \cite{KLMPP}: For any poset $P$, $\sat^*(n,P)$ is either bounded by a constant or grows at least as a logarithmic function of $n$. This was recently improved as follows.

\begin{theorem}[Freschi, Sharifzadeh, Spiga, Treglown \cite{FPST}]\label{dich}
For any poset $P$, either there exists a constant $c^*_P$ such that $\sat^*(n,P)\le c^*_P$ holds for all $n$ or we have $\sat^*(n,P)=\Omega(\sqrt{n})$.
\end{theorem}

In \cite{KLMPP}, it was shown that for any poset $P$ we have $\sat(n,P)\le c_P$ for some constant $c_P$ depending only on $P$ and not on $n$.

\subsection{Projective saturation}

The constant upper bounds of \cite{KLMPP} and Theorem \ref{dich} tell us that the ``(size of the) underlying set might not matter.'' This leads us to the idea to allow the saturated families $\cF$ to contain arbitrary sets, while keeping the aim that every $G\in \cB_n\setminus \cF$ should create a copy of the forbidden poset $P$. Will this change make strong saturation numbers constant, or will proofs become easier in this more general setting?

As $\sat(n,P)$ is always bounded by a constant,
from here on, all copies of any poset in this note are going to be strong copies unless otherwise stated. The formal definition is as follows. Note that $N$ does not really play any role.

\begin{defn}\label{def:pisat}
A family $\cF\subset \cB_N$ is \textit{projective $P$-saturated for $\cB_n$} if \begin{itemize}
    \item[(i)] 
    $\cF$ is $P$-free,
    \item[(ii)]
    for any $G\in \cB_n\setminus \cF$, the family $\cF\cup \{G\}$ contains a copy of $P$,
    \item[(iii)]
    for any two different $F, F'\in \cF$, we have $\pi(F):=F\cap [n] \neq F'\cap [n]=\pi(F')$.
\end{itemize}
The minimum size that a projective $P$-saturated family for $\cB_n$ can have is denoted by $\pisat(n,P)$. By definition, we have $\pisat(n,P)\le \sat^*(n,P)$.
\end{defn}

Observe that if we omit the last condition of Definition \ref{def:pisat}, then the problem becomes a lot less interesting.

\begin{proposition}\label{prop:sameproj}
    For any poset $P$, the saturation number satisfying only conditions (i) and (ii) of Definition \ref{def:pisat} is $|P|$ or $|P|-1$, depending on whether $P$ contains a smallest element or not, provided that $2^n$ is at least this big.
\end{proposition}


We will be interested in what statements on $\sat^*(n,P)$ remain true for $\pisat(n,P)$, and whether those that do, will need a new proof or not. Our first observation is that both the proof of the weaker dichotomy statement of \cite{KLMPP} and that of Theorem \ref{dich} work word-by-word in this scenario. 

\begin{theorem}\label{dich2}
For any poset $P$ either there exists a constant $c_{\pi,P}$ such that $\pisat^*(n,P)\le c_{\pi,P}$ holds for all $n$ or we have $\pisat^*(n,P)=\Omega(\sqrt{n})$.
\end{theorem}

The interested reader can find the proof of Theorem \ref{dich2} in
Appendix \ref{sec:app}.
On the other hand, we will show that the class of posets $P$ for which $\pisat(n,P)$ is constant is a \textit{strict} superclass of those for which $\sat^*(n,P)$ is constant.

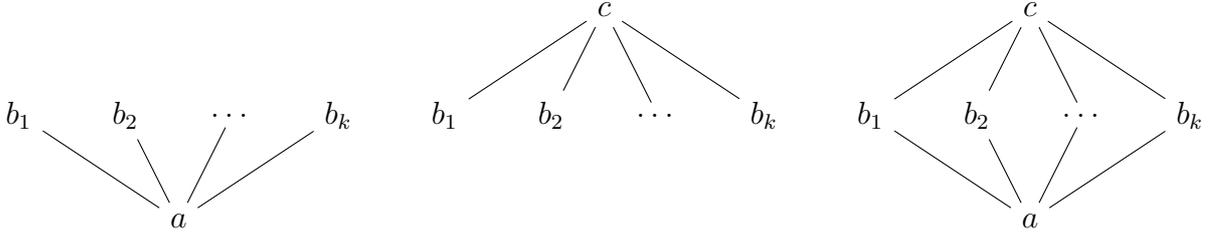
\begin{figure}
\begin{center}
    \begin{tikzpicture}[scale=0.7]
    \node (a) at (-15,0) {$b_1$};
  \node (b) at (-13,0) {$b_2$};
  \node (c) at (-11,0) {$\dots$};
  \node (d) at (-9,0) {$b_k$};
  \node (zero) at (-12,-2) {$a$};
  \draw (zero) -- (a); 
  \draw (b) -- (zero); 
  \draw (c) -- (zero); 
\draw (d) -- (zero);
\node (a) at (-7,0) {$b_1$};
  \node (b) at (-5,0) {$b_2$};
  \node (c) at (-3,0) {$\dots$};
  \node (d) at (-1,0) {$b_k$};
  \node (zero) at (-4,2) {$c$};
  \draw (zero) -- (a); 
  \draw (b) -- (zero); 
  \draw (c) -- (zero); 
\draw (d) -- (zero);
\node (a) at (1,0) {$b_1$};
  \node (b) at (3,0) {$b_2$};
  \node (c) at (5,0) {$\dots$};
  \node (d) at (7,0) {$b_k$};
  \node (zero) at (4,-2) {$a$};
  \node (zuro) at (4,2) {$c$};
  \draw (zero) -- (a); 
  \draw (b) -- (zero); 
  \draw (c) -- (zero); 
\draw (d) -- (zero);
\draw (zuro) -- (a); 
  \draw (b) -- (zuro); 
  \draw (c) -- (zuro); 
\draw (d) -- (zuro);
    \end{tikzpicture}
    \caption{Hasse diagrams of $\vee_k$, $\wedge_k$ and $D_k$}
    \label{fig:posets}
    \end{center}
\end{figure}

We denote by $A_k$ the $k$-element antichain, and by $C_k$ the $k$-element chain. The fork $\vee_k$, the cherry $\wedge_k$, and the diamond $D_k$ have point sets $\{a,b_1,b_2,\dots,b_k\}$, $\{b_1,b_2,\dots,b_k,c\}$, and $\{a,b_1,b_2,\dots,b_k,c\}$, respectively, with $a<b_i<c$ for all $i=1,2\dots,k$ and the $b_i$'s forming an antichain (see Figure \ref{fig:posets}).

\vskip 0.2truecm

All of $\wedge_k$, $\vee_k$, and $D_k$ have unbounded $\sat^*$-value if $k\ge 2$. This was shown in \cite{Fetal} using the so-called \textit{unique cover twin property} (UCTP). A poset $P$ has UCTP if whenever $q$ is the unique cover of $p$ in $z$, then there exists an element $z\in P$, not equal to $p$, of which the only cover is also $q$. The next proposition shows that $\pisat(n,P)$ can be bounded for posets with UCTP. For a poset $P$, we denote by $kP$ the poset obtained by $k$
pairwise incomparable copies of $P$. The unboundedness of $\sat^*(n,2C_2)$ (where $2C_2$ is the union of two incomparable chains on two elements each) was shown in \cite{KLMPP}.

\begin{proposition}\label{dia}\
\begin{enumerate}[(i)]
    \item 
    For any $k\ge 2$, we have $\pisat(n,\wedge_{k}),\pisat(n,D_{k})\le 2k+1$.
    \item
    We have $\pisat(n,2C_2)\le 8$
\end{enumerate}
\end{proposition}

On the other hand, we show that if $P=A_k$ for some $k\ge 2$, then there does not exist any projective $A_k$-saturated family that is not simply $A_k$-saturated.

\begin{proposition}\label{anti}
For any $k\ge 1$, if $\cF$ is projective $A_{k+1}$-saturated for $\cB_n$, then $\cF\subseteq \cB_n$. In particular, $\pisat(n,A_{k+1})=\sat^*(n,A_{k+1})$.
\end{proposition}

Although the UCTP proof of \cite{Fetal} does not remain valid in the context of $\pisat(n,P)$, with some extra work, however, one can save the argument for $\vee_k$. The next proposition, together with Proposition \ref{dia}, shows that another phenomenon of $\sat^*$ does not hold for $\pisat$: if $P^D$ is the dual of $P$ obtained by reversing all relations of $P$, then by taking complements of $P$-saturated families,  it is clear that $\sat^*(n,P)=\sat^*(n,P^D)$ holds for all $P$, but we can have $\pisat^*(n,P)\ne\pisat^*(n,P^D)$.

\begin{proposition}\label{vee}
For any $n$, we have $\pisat(n,\vee)=\sat^*(n,\vee)=n+1$.

For any $k\ge 2$, 
we have $\pisat(n,\vee_k)=\Omega(\sqrt{n})$.
\end{proposition}

\subsection{External saturation}

As $\pisat(n,\vee)\ne \pisat(n,\wedge)$, it is natural to attempt to find a more symmetric variant of projective saturation, for example, by fixing some $A\cB_n=\{A\cup B : B\in \cB_n\}$ for some $A\subset \{n+1,\ldots,N\}$,
and requiring saturation only in $A\cB_n$.
Again, $N$ does not really play any role.

\begin{defn}\label{def:extsat}
A family $\cF\subset \cB_N$ is \textit{external $P$-saturated for $\cB_n$} if there exists $A\subseteq [n+1,N]$ such that
\begin{itemize}
    \item[(i)] 
    $\cF$ is $P$-free,
    \item[(ii)]
    for any $G\in A\cB_n\setminus \cF$, the family $\cF\cup \{G\}$ contains a copy of $P$,
    \item[(iii)]
    for any two different $F, F'\in \cF$, we have $\pi(F):=F\cap [n] \neq F'\cap [n]=\pi(F')$.
\end{itemize}
The minimum size that an external $P$-saturated family for $A\cB_n$ can have is denoted by $\extsat(n,P)$. Observe that if $\cF\subseteq \cB_N$ is external $P$-saturated for $A\cB_n$, then $\overline{\cF}=\{[N]\setminus F:F\in\cF\}$ is external $P^D$-saturated for $A'\cB_n$, where $A'=[N]\setminus ([n]\cup A)$. Thus, we have $\extsat(n,P)\le \min\{ \pisat(n,P),\pisat(n,P^D)\}$.
\end{defn}

The proof of the dichotomy results Theorem \ref{dich} and Theorem \ref{dich2} remain valid for external saturation, but since surprisingly we do not have an example of $P$ with unbounded external saturation number, we have the following conjecture.

\begin{conjecture}\label{conj:extsat}
    $\extsat(n,P)\le C_P$ for every $P$ for some $C_P$ independent of $n$.
\end{conjecture}
    

However, it is easy to see that the class of posets for which $\extsat(P,n)$ is bounded by a constant is a strict superclass of those posets for which $\pisat(n,P)$ is a constant. This follows by considering $\vee_k$ for which $\pisat(n,\vee_k)=\Omega(\sqrt{n})$ by Proposition \ref{vee}, but $\extsat(n,\vee_k)\le \pisat(n,\wedge_k)\le 2k+1$ by Proposition \ref{dia}. It is again not hard to see that there exist posets $P$ for which both $\pisat(n,P)$ and $\pisat(n,P^D)$ are unbounded, but $\extsat(n,P)$ is bounded by a constant. The simplest such poset is the antichain $A_k$ for any $k\ge 2$.

\begin{proposition}
    $\extsat(n,A_k)=k-1$ for all $k$ and $n$.
\end{proposition}
\begin{proof}
    Let $N=n+2$ and $A=\{n+1\}$. Take any antichain $\mathcal{X}$ with $k-1$ elements from $\cB_n$, and let $\cF=\{F\cup \{n+2\}\mid F\in \mathcal{X}\}$. This is clearly $A_k$-free, and it is also $A_k$-saturated for $A\cB_n$, as each element of $A\cB_n$ is incomparable to each element of $\cF$.
\end{proof}

Our final contribution is a class of height 2 posets for which we can prove Conjecture \ref{conj:extsat} by constructing external saturated families of constant size. It is easy to see that if the comparability graph of $P$ contains an isolated vertex or an isolated edge, then $\extsat(n,P)= |P|-1$ provided $2^n\ge |P|-1$ (we do not need $h(P)=2$ for this). Indeed, we can have $N=n+2$, $A=\{n+1\}$, then if $p\in P$ is an isolated vertex and $\cG\subseteq 2^{[n]}$ is a copy of $P\setminus \{p\}$, then $\{G\cup \{n+2\}: G\in \cG\}$ is external $P$-saturated for $A\cB_n$. Similarly, if $P$ has a $C_2$-component, and $\cG\subseteq 2^{[n]}$ is a copy of $P\setminus C_2$, then $\{\{n+1\}\}\cup \{G\cup \{n+2\}:G\in \cG\}$ is external $P$-saturated for $A\cB_n$. 

Let $P=P_1\cup P_2$ be a height two poset with $P_1$ being the set of minimal elements of $P$ and $P_2=P\setminus P_1$ being the set of maximal, but not also minimal, elements. We say that a height 2 family $\cF=\cF_1\cup \cF_2 \subset \cB_n$ is \textit{almost saturated} if it is $P$-free and adding any $F\in \cB_n \setminus \cF$ creates  a weak copy of $P$ that is almost strong, in the following sense.
If $F$ is on the top level of the copy, it might have extra containment relations with the members of $\cF_2$, while if $F$ is on the bottom level, it might have extra containment relations with the members of $\cF_1$, but otherwise all containments are exactly the same as in $P$.

\begin{proposition}\label{prop:extsat}
Suppose that for some height two, isolated $C_2$-free poset $P$ we have an almost saturated family $\cF=\cF_1\cup\cF_2$.
Then $\extsat(n,P)\le 2|\cF|$ for every $n$.
\end{proposition}

This can be applied to $K_{s,t}$, the poset whose Hasse-diagram is a complete bipartite graph on $s+t$ vertices.

\begin{corollary}\label{cor:Kst}            $\extsat(n,K_{s,t})\le 4(s+t-1)$ for all $s,t,n$.
\end{corollary}

We could also prove Conjecture \ref{conj:extsat} for some other specific posets, but we do not have any other interesting general result.

\section{Proofs}

\begin{proof}[Proof of Proposition \ref{prop:sameproj}]
    If $2^n<|P|$, then $P\nsubseteq \cB_n$, so $\cF= \cB_n$ is saturated.

    Otherwise, $|P|-1$ is clearly a lower bound, as if $\cF$ consists of at most $|P|-2$ sets, then $\cF\cup \{G\}$ cannot contain any copies of a poset of size $|P|$. Also, if $P$ does not contain a smallest element, then the empty set does not help to create any copy of $P$, and thus it must be an additional element of any $P$-saturated family.

    On the other hand, fix a minimal element $p\in P$ and a copy $\cF\subset \cB_N$ of $P\setminus \{p\}$ such that $\emptyset\ne F\subset \{n+1,\ldots,N\}$ for all ${F\in \cF}$, i.e., $\pi(F)=\emptyset$ for all $F\in \cF$ but $\emptyset\notin \cF$.
    Also fix a bijection $f:P\setminus \{p\}\rightarrow \cF$.
    Then $\cG:=\{f(q):p\not\leqslant q\}\cup \{[n]\cup f(q):p\leqslant q\}\cup \{\emptyset\}$ $P$-saturates any $H\in \cB_n\setminus \cG$ with $H$ playing the role of $p$ for any $H$.
    For example, if $P=\wedge$, then we can take $\cF=\{\{n+1\},\{n+1,n+2\}\}$, and
    $\cG=\{\{n+1\},[n+2],\emptyset\}$.
    Moreover, if $P$ has a smallest element, then even $\cG\setminus\{\emptyset\}$ will $P$-saturate --- note that in this case all elements of $\cG$ contain $[n]$.
\end{proof}

\begin{proof}[Proof of Proposition \ref{dia}]
The same construction works for $D_k$ and $\wedge_k$, proving (1). Let $$\cF_{n,k}:=\left\{\emptyset,[n+k-1],[2k-2]\cup\{n+k\}\right\}\cup\left\{\{i,n+i\},\{k-1+i,n+k\}: i=1,2,\dots k-1\right\}.$$
Observe that $\cF_{n,k}$ consists of 2 incomparable $\wedge_{k-1}$ posets, plus the $\emptyset$, so it is $\wedge_k$- and thus $D_k$-free.
Furthermore, if $F\in \cB_n\setminus (\{\emptyset\}\cup \{\{i\}:i=1,2,\dots,k-1\})$, then $F,[n+k-1],\{\{i,n+i\}:i=1,2,\dots,k-1\}$ form a copy of $\wedge_k$ (and together with $\emptyset$ a copy of $D_k$). Finally, $\{j\},[2k-2]\cup \{n+k\}, \{\{k-1+i,n+k\}: i=1,2,\dots,k-1\}$ form a copy of $\wedge_k$ (together with $\emptyset$ a copy of $D_k$) for any $j=1,2,\dots,k-1$. So $\cF_{n,k}$ is indeed projective $\wedge_k$- and $D_k$-saturated for $\cB_n$.\\

To prove (2), we can take the followng projective $2C_2$-saturated family of size 8:
$$\emptyset, \{3\}, \{1, 2\}, \{2, 3\}, \{2,n+1\}, \{1, 2, 3, n+1\},  [n]\setminus\{2\}, [n+1]\setminus\{1\}.$$
It can be verified with a similar simple case analysis, which we omit here, that this family is indeed projective $2C_2$-saturated for $\cB_n$.
\end{proof}

\begin{proof}[Proof of Proposition \ref{anti}]
Suppose that there exists a family $\cF$ that is projective $A_{k+1}$-saturated for $\cB_n$ such that $F\in \cF\setminus \cB_n$. Let us pick $F$ such that $|F|$ is minimal over all $F'\in \cF\setminus \cB_n$. As $\cF$ is $A_{k+1}$-free, one can partition $\cF$ into $k$ many chains $\cC_1=\{C^1_1\subset C^1_2\subset \dots \subset C^1_{i_j}\},\cC_2,\dots, \cC_k$, and assume that $F$ is the $j$th set in $\cC_1$. Then by the minimality of $|F|$, we have $C^1_{j-1}\in \cB_n$, and as $\pi(A)\neq \pi(B)$ for all $A,B\in \cF$, we must have $C^1_{j-1}\subsetneq \pi(F)$. Also, $\pi(F)\notin \cF$ as $F\in \cF\setminus \cB_n$. But then $\pi(F)$ can be added to $\cC_1$ and so $\cF\cup \{\pi(F)\}$ is still a union of $k$ chains and thus $A_{k+1}$-free. This contradiction finishes the proof.
\end{proof}

\begin{proof}[Proof of Proposition \ref{vee}] 
For the upper bound, observe that for any $k\le n$ we have $\pisat(n,\vee_k)\le \sat^*(n,\vee_k)\le n+1$ shown by the family $\binom{[n]}{\ge n-1}=\{F\subseteq [n]:|F|\ge n-1\}$. 

The beginning of the proof of both lower bounds is the same.
Let $\cF$ be a projective $\vee_k$-saturated family for $\cB_n$. We write $\cF^{in}=\cF\cap \cB_n$ and $\cF^{out}=\cF\setminus \cB_n$. We start with a claim that states that every pair $A\in \cF^{in}, B\in \cF^{out}$ is incomparable. Clearly, we cannot have $B\subseteq A$.

\begin{claim}\label{inandout}
There does not exist any pair $A\in \cF^{in}, B\in \cF^{out}$ with $A\subseteq B$.
\end{claim}

\begin{proof}
Suppose there are such $A$ and $B$, and consider such a pair that minimizes $|\pi(B)|$ over all such pairs. As $B\in \cF$, we have $\pi(B)\notin \cF$, so there exist $F_1,F_2,\dots,F_k\in \cF$ such that $\pi(B),F_1,\dots,F_k$ form a copy of $\vee_k$. If $\pi(B)$ is the root of this copy, then $A,F_1,\dots,F_k$ form a copy of $\vee_k$ in $\cF$---a contradiction. If $\pi(B)$ is not the root but, say, $F_1$ is, then the only way how $B,F_1,F_2,\dots,F_k$ would not form a copy of $\vee_k$ is $F_i\subset B$ for some $2\le i\le k$. But then $F_i\notin \cB_n$ would contradict the choice of $B$, while $F_i\in \cB_n$ would mean $F_i\subset \pi(B)$ contradicting that $F_1,F_2,\dots,F_k,\pi(B)$ form a copy of $\vee_k$ with $F_1$ being the root, as then $F_i$ and $\pi(B)$ should be incomparable.
\end{proof}

Claim \ref{inandout} implies that if $\emptyset \in \cF$, then $\cF=\cF^{in}$ and thus $|\cF|\ge \sat^*(n,\vee_k)$.

\begin{claim}\label{n}
The set $[n]$ belongs to $\cF$.
\end{claim}

\begin{proof}
Suppose not. Then $[n],F_1,F_2,\dots,F_k$ is a copy of $\vee_k$ for some $F_1,F_2,\dots,F_k\in \cF$. If $[n]$ is the root, then $\pi(F_i)=[n]$ for all $i$---a contradiction as $k\ge 2$. If, say, $F_1$ is the root, then $F_1\in \cF^{in}$ as $F_1\subseteq [n]$, while for all other $i$ we must have $F_i\in\cF^{out}$ as $[n]$ and $F_i$ must be incomparable. But then $F\subset F_i$ contradicts Claim \ref{inandout}.
\end{proof}

Next, we prove $\pisat(n,\vee_k)=\Omega(\sqrt{n})$.
As in \cite{Fetal}, it is enough to prove that for any $x,y\in [n]$, there exists $F\in \cF$ with $|F\cap \{x,y\}|=1$. We borrow ideas from \cite[Lemma 8]{Fetal}. Suppose towards a contradiction that for some $x,y\in [n]$, we have that for any $F\in \cF$ the size of $F\cap \{x,y\}$ is either 0 or 2. We will show that then $F\cap \{x,y\}=\emptyset$ for all $F\in\cF$, which would contradict Claim \ref{n}. So let us fix a minimum size $S$ such that $S_{x,y}:=S\cup \{x,y\}\in\cF$. 

Observe first that whenever $F\in \cF$ is incomparable to $S_x:=S\cup \{x\}$, then so are $F$ and $S_{x,y}$. Indeed, if $|F\cap \{x,y\}|=0$, then the incomparability of $F$ and $S_x$ yields the existence of an element $z\in F\setminus S_{x,y}$ and thus $F$ and $S_{x,y}$ are incomparable. If $F\supset \{x,y\}$, then by choice of $S$, we have $|F|\ge |S_{x,y}|$, and thus the only possibility for comparability is $F\supset S_{x,y}\supset S_x$, which contradicts the incomparability of $F$ and $S_x$.

As $F\cap \{x,y\}$ has size 0 or 2 for all $F\in \cF$, we know that $S_x\notin \cF$, and thus there exists $\cF'\subseteq \cF$ such that $\{S_x\}\cup \cF'$ is a copy of $\vee_k$. We claim that if $F\in \cF'$ covers $S_x$ in the poset structure of $\{S_x\}\cup \cF'$, then $F$ covers $S_{x,y}$ in the poset structure of $\{S_{x,y}\}\cup \cF'$. Indeed, the only problem could be when $S_x$ is the root of the copy $\{S_x\}\cup \cF'$ of $\vee_k$. So all $F\in \cF'$ contains $x$ and thus $y$. As they are incomparable, they all must properly contain $S_{x,y}$. These statements and the previous paragraph imply that $\{S_{x,y}\}\cup\cF'$ is a copy of $\vee_k$, which contradicts the $\vee_k$-free property of $\cF$.\\

Next, we prove the lower bound $\pisat(n,\vee)\ge n+1$. Consider the family $\pi(\cF^{out})=\{\pi(F):F\in \cF^{out}\}$ and the downset $\cD$ generated by $\pi(\cF^{out})$, i.e., $\cD=\{G: \exists F\in \cF^{out} ~\text{such that}\ G\subseteq \pi(F)\}$. Observe that $\cF\cap \cD=\emptyset$ by Claim \ref{inandout}. 
So $\cF^{in}\subset \cU:=\cB_n\setminus \cD$. Also, for any $G\in \cU\setminus \cF$, in the copy $G,F,F'$ of $\vee$, we must have $F,F'\in \cF^{in}$. Indeed, by the definition of $\cU$, $G$ is only comparable to sets in $\cF^{in}$, and Claim \ref{inandout} ensures that then both $F,F'$ must belong to $\cF^{in}$.\\

The remainder of the proof is a slight modification of that in \cite{Fetal} showing $\sat^*(n,\vee)=n+1$, which we include here for completeness. Our goal is to define an injection $f:\binom{[n]}{n-1}\rightarrow \cF\setminus \{[n]\}$. If we manage to do so, then Claim \ref{n} ensures $|\cF|\ge n+1$. If for some $G\in \binom{[n]}{n-1}$ there exists $F\in \cF$ with $\pi(F)=G$ (in particular, if $G\in \cF$), then we let $f(G)=F$. Otherwise, 
there exist $F,F'\in \cF^{in}$ that form a copy of $\vee$ with $G$. As $[n]$ is the only one set in $\cF^{in}$ that contains $G$, $G$ cannot be the root of this copy, so we can assume that the root is $F$. Let us pick $F,F'$ such that $|F'|$ is minimal. 
We let $f(G)=F$.

To show that $f$ is an injection, denote by $x$ the element for which $G=[n]\setminus \{x\}$.
Following \cite{Fetal}, we show that $F'=F\cup\{x\}$.
As $F'\nsubseteq G$, we have $x\in F'$.
Define $F_{\text{-}x}'=F'\setminus\{x\}$.
Suppose for contradiction that $F\neq F_{\text{-}x}'$.
It is clear that $F\subseteq F_{\text{-}x}'$, as $x\notin F\subsetneq F'$.
By the minimality of $F'$, then $F_{\text{-}x}'\notin \cF$.
Since $\cF$ is projective $\vee$-saturated, $\cF\cup \{F_{\text{-}x}'\}$ contains a copy of $\vee$.
As $F\in \cF$ and $F\subsetneq F_{\text{-}x}'$, in this copy $F_{\text{-}x}'$ cannot be the smallest element, or we could replace it with $F$, contradicting that $\cF$ is $\vee$-free.
So we have some $H\subsetneq F_{\text{-}x}'$ and $H'\nsubseteq F_{\text{-}x}'$ in $\cF$ such that $H\subsetneq H'$.
As $F',H,H'$ cannot form a $\vee$, we have $H'\subsetneq F'$.
Since $H'\nsubseteq F_{\text{-}x}'$ but $H'\subset F'=F_{\text{-}x}'\cup\{x\}$, we can conclude that $x\in H'$ and thus $H'\nsubseteq G$.
But $H\subset F_{\text{-}x}'\subset G$, so $G,H,H'$ form a $\vee$, contradicting the minimality of $|F'|$.
Thus, $F'=F\cup\{x\}$.
In other words, if $F=f([n]\setminus \{x\})$ for some $x$, then $F\cup\{x\}\in \cF$.

To finish the proof of the injectivity of $f$, assume for a contradiction that $f([n]\setminus \{x\})=F=f([n]\setminus \{y\})$.
But then $F,F\cup\{x\},F\cup\{y\}$ form a $\vee$, contradiction.
\end{proof}

\begin{proof}[Proof of Proposition \ref{prop:extsat}]
    Let $N=n+3$, and set $A=\{n+3\}$.
    Our external $P$-saturated family $\cG$ will consist of four parts, $\cG=\cG_1\cup\cG_1'\cup\cG_2'\cup\cG_2$, defined as follows, using the shorthand $\cF+A=\{F\cup A\mid F\in \cF\}$.\\ 
    $\cG_1=\cF_1$,\\
    $\cG_1'=\cF_1+\{n+1\}$,\\
    $\cG_2'=\cF_2+\{n+1,n+2\}$,\\
    $\cG_2=\cF_2+\{n+1,n+2,n+3\}$.\\
    Note that $\cG\cap A\cB_n=\emptyset$, that the relations among $\cG_1$ and $\cG_1'$ form $|\cF_1|$ isolated $C_2$'s, while the relations among $\cG_2$ and $\cG_2'$ form $|\cF_2|$ isolated $C_2$'s, and the relations among $\cG^*_1$ and $\cG^*_2$ are the same as among $\cF_1$ and $\cF_2$, where $\cG^*_i$ is either $\cG_i$ or $\cG'_i$.

    First, we show that $\cG$ is $P$-free.
    Suppose for contradiction that it contains a $g(P)$ copy of $P$.
    If for some $F\in \cG_1\cap g(P)$ we have $F\cup\{n+1\}\in \cG_1'\cap g(P)$, then these would form an isolated $C_2$ in $g(P)$, as they are comparable to the same sets of $\cG_2$ and $\cG_2'$, and $P$ has only two levels.
    But this is impossible as $P$ is $C_2$-free.
    Otherwise, we can suppose that $g(P)\cap \cG_i'=\emptyset$, as we can replace each set with a respective one from $\cG_i$.
    Thus, $g(P)\subset \cG_1\cup \cG_2$, but $\cG_1$ and $\cG_2$ have the same relations as $\cF_1$ and $\cF_2$, and $\cF$ is $P$-free.
    Therefore, $g(P)$ cannot be a copy of $P$.
    
    Next, suppose that we add some $G=F\cup A\in A\cB_n$ to $\cG$.
    Consider an $f(P)$ weak but almost strong copy of $P$ that $F\in \cB_n$ creates when we add it to $\cF$.
    We will show how to turn this into a $g(P)$ strong copy of $P$ for $\cG$.
    Without loss of generality, assume that $F$ is on the top level of $f(P)$.
    Then $G$ will also be on the top level.
    For all other top level sets in $f(P)$, add $\{n+1,n+2\}$ to them, obtaining sets from $\cG_2'$.
    The bottom level of $f(P)$ remains unchanged, those sets are from $\cG_1$.
    This finishes the description of $g(P)$.

    As the containment relations among $\cG_1$ and $\cG_2'$ are the same as the relations among $\cF_1$ and $\cF_2$.
    The contaiment relations of $G$ and $\cG_1$ are also the same as the relations among $F$ and $\cF_1$.
    Finally, there are no containments between $G$ and $\cG_2$.
    Thus, $g(P)$ is indeed a strong copy of $P$, as $f(P)$ was an almost strong copy.  
\end{proof}

\begin{proof}[Proof of Corollary \ref{cor:Kst}]
    If $s=t=1$, the statement is trivial.\\
    Othwerwise, we can choose $\cF=\cF_1\cup\cF_2$, where\\
    $\cF_1=\{
    \{1\},\ldots,\{s+t-1\}\}$,\\
    $\cF_2=\{
    [n]\setminus\{1\},\ldots,[n]\setminus\{s+t-1\}\}$.\footnote{This is not the best family for weak saturation; in \cite{KLMPP} it was shown that $\sat(n,K_{s,t})\le 2(s+t-1)-2$.}\\
    
    By the pigeonhole principle, $\cF_1\cup\cF_2$ is $K_{s,t}$-free.
    
    Adding any $F\in \cB_n$ to $\cF$ will create an almost strong copy of $K_{s,t}$.
    Indeed, if $|F\cap [s+t-1]|\ge s$, then $F$ contains some $s$ sets of $\cF_1$, which are further contained in some $t-1$ sets of $\cF_2$,
    while if $|F\cap [s+t-1]|\le s-1$, then $F$ is contained in some $t$ sets of $\cF_2$, which further contain some $s-1$ sets of $\cF_1$.
\end{proof}    
    

\appendix
\section{}\label{sec:app}

\begin{proof}[Proof of Theorem \ref{dich2}] As we mentioned, the proof follows word-by-word that in \cite{FPST}. Suppose first that for  every $i\in n$ and every $\cF$ that is projective $P$-saturated for $\cB_n$, there exist $F,F'\in \cF$ with $F\setminus F'=\{i\}$. Then, clearly, $|\cF|(|\cF|-1)\ge n$ and thus $|\cF|>\sqrt{n}$.

So we can assume that there exists $n_0$ and a family $\cF$ that is projective $P$-saturated for $\cB_{n_0}$ such that $\{n_0\}\neq F\setminus F'$ for any $F,F'\in \cF$. Although in the definition of projective saturation, we required $\cF\subseteq \cB_N$ for some arbitrary large $N$, for convenience we assume that $\bigcup_{F\in \cF}F=[n_0]\cup A$, where $A$ is disjoint with the set of integers. Then for any $n\ge n_0$, we define a bijection $f=f_n: 2^{A\cup [n_0]}\rightarrow 2^{A \cup [n]}$ that ``blows up'' $n_0$ into $n-n_0$ elements as
\[
f(F)=\begin{cases}
    F & \text{if $n_0\notin F$, } \\
    F \cup [n_0+1,n] & \text{if $n_0\in F$, }
\end{cases}
\]
and let $\cF_n=\{f(F):F \in \cF\}$.
Clearly, $|\cF_n|=|\cF|$. So if we can prove that $\cF_n$ is projective $P$-saturated for $\cB_n$, then $\pisat(n,P)\le \max\{\pisat(m,P):m\le n_0\}$. By definition, the poset structure of $\cF$ and $\cF_n$ are the same, so $\cF_n$ is $P$-free as so is $\cF$.

We are left to prove that for any $G\in 2^{[n]}\setminus \cF_n$ there exists a copy of $P$ in $\cF\cup \{G\}$. We distinguish four cases.

\medskip

\textsc{Case I.} $G \subseteq [n_0-1]$

\smallskip

Then $f(G)=G$, and so $G\notin \cF$ as otherwise it would also belong to $\cF_n$. So $\cF\cup \{G\}$ contains a copy $\cG$ of $P$ and as $f$ is a poset-isomorphism, $f[\cG]$ is a copy of $P$ in $\cF_n \cup \{G\}$.

\medskip

\textsc{Case II.} $G\supseteq [n_0,n]$

\smallskip

Then for $G^*:=G\cap [n_0]\subset [n_0]$, we have $f(G^*)=G$ and thus $G^*\notin \cF$, and so $\cF\cup \{G^*\}$ contains a copy $\cG$ of $P$. But then $f[\cG]$ is a copy of $P$ in $\cF_n\cup \{G\}$.

\medskip

\textsc{Case III.} $n_0\in G, \exists j\in [n_0+1,n]\setminus G$

\smallskip

We need the following claim.

\begin{claim}\label{klem}
    Either $\cF_n\cup \{G\}$ contains a copy of $P$ or

    (1) there exists $F_1\in \cF$ with $n_0\in F_1\subseteq G$,

    and

    (2) there exists $F_2\in \cF$ with $n_0\notin F_2$, $G\cap [n_0-1]\subseteq F_2$.
\end{claim}

\begin{proof}
    Suppose first we do not have $F_1\in \cF$ satisfying (1). In particular, $G_1:=G\cap [n_0]\notin \cF$, so $\cF\cup \{G_1\}$ contains a copy $\cG_1$ of $P$ with $G_1\in \cG_1$. Now the containment relation of $G_1$ and any $F\in\cF$ is the same as that of $G$ and $f(F)$ unless $F$ is a subset of $G_1\subset G$ containing $n_0$. So either we find $F_1$ with the required property, or $\{G\}\cup f[\cG_1\setminus \{G_1\}]$ is a copy of $P$ in $\cF_n \cup \{G\}$. 
    
    Suppose next we do not have $F_2\in \cF$ satisfying (2).  In particular, $G_2:=G\cap [n_0-1]\notin \cF$, so $\cF\cup \{G_2\}$ contains a copy $\cG_2$ of $P$ with $G_2\in \cG_2$.  Now the containment relation of $G_2$ and any $F\in\cF$ is the same as that of $G$ and $f(F)$ unless $F$ is a superset of $G_2$ not containing $n_0$. So either we find $F_2$ with the required property or $\{G\}\cup f[\cG_2\setminus \{G_2\}]$ is a copy of $P$ in $\cF_n \cup \{G\}$. 
\end{proof}

As Claim \ref{klem} implies either the existence of $F_1,F_2$ with $F_1\setminus F_2=\{n_0\}$ or a copy of $P$ in $\cF_n\cup \{G\}$, we conclude to the latter as the assumption of Theorem \ref{dich2} states that there are no such $F_1,F_2$ in $\cF$.

\medskip

\textsc{Case IV.} $n_0\notin G, \exists j\in [n_0+1,n]\cap G$

\smallskip

Consider $G^*:=G\setminus \{j\} \cup \{n_0\}$. By \textsc{Case III}, $\cF_n \cup \{G^*\}$ contains a copy $\cG^*$ of $P$ with $G\in \cG^*$. All sets in $\cF_n$ are either disjoint to $[n_0,n]$, or contain $[n_0,n]$, so $G$ and $G^*$ have the same comparability to all sets in $\cF_n$, therefore $\cG^*\setminus \{G^*\}\cup \{G\}$ is a copy of $P$ in $\cF_n\cup \{G\}$.
\end{proof}

\begin{thebibliography}{99}
\bibitem{Betal}
Bastide, P., Groenland, C., Jacob, H., Johnston, T. (2022). Exact antichain saturation numbers via a generalisation of a result of Lehman-Ron. arXiv preprint arXiv:2207.07391.
\bibitem{DI}
Đanković, I., Ivan, M. R. (2023). Saturation for Small Antichains.
The Electronic Journal of Combinatorics 30(1), P1-3.
\bibitem{Fetal}
Ferrara, M., Kay, B., Kramer, L., Martin, R. R., Reiniger, B., Smith, H. C., Sullivan, E. (2017) The saturation number of induced subposets of the Boolean lattice. Discrete Mathematics, 340(10), 2479-2487.
\bibitem{FPST}
Freschi, A., Piga, S., Sharifzadeh, M., Treglown, A. (2022). The induced saturation problem for posets. arXiv preprint arXiv:2207.03974.
\bibitem{Getal}
Gerbner, D., Keszegh, B., Lemons, N., Palmer, C., Pálvölgyi, D., Patkós, B., Saturating Sperner families,
Graphs and Combinatorics 29 (5), 1355-1364
\bibitem{GP}
Gerbner, D., Patkós, B. (2018). \textit{Extremal finite set theory}. Chapman and Hall/CRC.
\bibitem{GL}
Griggs, J. R., Li, W. T. (2016). Progress on poset-free families of subsets. In Recent trends in combinatorics, 317--338. Springer, Cham.
\bibitem{I}
Ivan, M. R. (2020). Saturation for the butterfly poset. Mathematika, 66(3), 806-817.
\bibitem{I2}
Ivan, M. R. (2022). Minimal Diamond-Saturated Families. Contemporary Mathematics, 3(2), 81.
\bibitem{KLMPP}
Keszegh, B., Lemons, N., Martin, R. R., Pálvölgyi, D., Patkós, B. (2021). Induced and non-induced poset saturation problems. Journal of Combinatorial Theory, Series A, 184, 105497.
\bibitem{MSW}
Martin, R. R., Smith, H. C.,  Walker, S. (2019). Improved bounds for induced poset saturation. The Electronic Journal of Combinatorics 27(2), P2-31.
\bibitem{MNS}
Morrison, N., Noel, J. A., Scott, A. (2014). On Saturated $ k $-Sperner Systems. The Electronic Journal of Combinatorics 21(3), P3-22.
\end{thebibliography}
\end{document}